\documentclass[11pt]{article}
\usepackage{bbm}
\usepackage{mathrsfs}
\usepackage{amsmath}
\usepackage{amsthm}
\usepackage{amsfonts}
\usepackage{amssymb}
\usepackage{latexsym}
\usepackage{amstext}
\usepackage{array}
\usepackage{graphicx}
\usepackage{MnSymbol}
\usepackage[style=plaintop]{floatrow}

\newtheorem{thm}[subsection]{Theorem}
\newtheorem{cor}[subsection]{Corollary }
\newtheorem{Def}[subsection]{Definition}
\newtheorem{lem}[subsection]{Lemma}
\newtheorem{remark}[subsection]{Remark}

\newtheorem{prop}[subsection]{Proposition}
\newtheorem{exm}[subsection]{Example}
\newcommand{\bthm}{\begin{thm} }
\newcommand{\ethm}{\end{thm} }
\newcommand{\bpro}{\begin{prop}}
\newcommand{\epro}{\end{prop}}
\newcommand{\bdf}{\begin{Def}}
\newcommand{\edf}{\end{Def}}
\newcommand{\bexm}{\begin{exm}}
\newcommand{\eexm}{\end{exm}}
\newcommand{\blem}{\begin{lem}}
\newcommand{\elem}{\end{lem}}
\newcommand{\bpf}{\begin{proof}}
\newcommand{\epf}{\end{proof}}
\newcommand{\bcor}{\begin{cor}}
\newcommand{\ecor}{\end{cor}}
\newcommand{\beq}{\begin{equation}}
\newcommand{\eeq}{\end{equation}}
\newcommand{\ba}{\begin{array}}
\newcommand{\ea}{\end{array}}
\newcommand{\bea}{\begin{eqnarray}}
\newcommand{\eea}{\end{eqnarray}}

\newcommand{\brem}{\begin{remark}}
\newcommand{\erem}{\end{remark}}

\begin{document}

%\footnote[0] { }

\title{Primality test for numbers of the form $(2p)^{2^n}+1$ }
\author{
Yingpu Deng and Dandan Huang\\\\
Key Laboratory of Mathematics Mechanization,\\
NCMIS, Academy of Mathematics and Systems Science,\\
Chinese Academy of Sciences, Beijing 100190, P.R. China\\
Email: \{dengyp, hdd\}@amss.ac.cn}

\date{}
\maketitle
\begin{abstract}
We describe a primality test for number $M=(2p)^{2^n}+1$ with odd
prime $p$ and positive integer $n$. And we also give the special primality
criteria for all odd primes $p$ not exceeding $19$. All these
 primality tests run in polynomial time in log$_{2}(M)$. A certain special $2p$-th
reciprocity law is used to deduce our result.

\end{abstract}

\section{Introduction}

Primality testing is an important problem in computational number theory. Although this has been proved to be a \textbf{P} problem by Agrawal, Kayal and Saxena \cite{AKS} in 2004, finding more efficient algorithms for specific families of numbers makes yet a lot of sense. Let $a>1$ be an integer. Considering the special family of prime numbers of the form $a^n\pm1$, when $a$ is fixed and $n>0$ is varied, is then a natural problem. For $a^n-1$, it is easy to see that it suffices to consider the case when $a=2$ and $n=p$ is a prime. Numbers of the form $2^p-1$ are called Mersenne numbers. For these numbers, there is a famous primality test named Lucas-Lehmer test given by Lucas \cite{lu} and Lehmer \cite{le}. Here, we recall it:

    \textbf{Lucas-Lehmer test.}\quad Let $M_p=2^p-1$ be Mersenne number, where $p$ is an odd prime. Define a sequence $\{u_k\}$ as follows: $u_0=4$ and $u_k=u_{k-1}^2-2$ for $k\geq1$. Then $M_p$ is a prime if and only if $u_{p-2}\equiv0\pmod{M_p}$.

For $a^n+1$, it is easy to see that it suffices to consider the case when $a$ is even and $n$ is a power of 2. When $a=2$, numbers of the form $2^{2^n}+1$ are called Fermat numbers. For these numbers, there is also a primality test due to P\'{e}pin (see \cite{Williams}):

\textbf{P\'{e}pin test.}\quad Let $F_n=2^{2^n}+1$ be the $n$-th Fermat number, with $n>0$. Then $F_n$ is a prime if and only if $3^{(F_n-1)/2}\equiv-1\pmod{F_n}$.

In this paper, we consider the primality of $M=(2p)^{2^n}+1$, where $p$ is an odd prime. For $p=3$ and $p=5$, Williams obtained primality tests for them using Lucas functions in \cite{wi}. However, for a general $p$, it seems that there is no any known work to deal with the primality of these numbers. Notice that, from the work in \cite{dl}, it is possible to give primality tests for these numbers. However, quoting directly the work in \cite{dl}, it will give primality tests whose seeds of sequences will depend on $M$, and this is not our desire. We can do better in this paper, that is, using a certain special $2p$-th reciprocity law, we can give primality tests whose seeds of sequences will not depend on $M$, at least for $p\leq19$.

This paper is organized as follows. In Section \ref{SecPre} we give the definition of power residue symbol and prove a certain special $2p$-th reciprocity law that will be used in later section. In Section \ref{SecMain} we state and prove our main result. In Section \ref{SecExpli}, we give explicit primality tests for $M=(2p)^{2^n}+1$ for each odd prime $p\leq19$. In Section \ref{SecImple} we give the implementation and computational results for $p=3,5$.

%\smallskip

\section{Preliminaries}\label{SecPre}

What we state in this section may be found in \cite[Chapter 14]{ir}.

For a positive integer $m$, let $\zeta_m=e^{2\pi\sqrt{-1}/m}$ be the complex primitive $m$-th root of unity, and
$D=\mathbb{Z}[\zeta_m]$ the ring of integers of the cyclotomic field $\mathbb{Q}(\zeta_m)$. Let $\mathfrak{p}$ be a prime ideal of $D$ lying over a
rational prime $p$ with gcd$(p, m)=1$. For every $\alpha\in D$, the $m$-th power residue
symbol $\left(\frac{\alpha}{\mathfrak{p}}\right)_m$ is defined by:

$\mathrm{(1)}$ If $\alpha\in\mathfrak{p}$, then $\left(\frac{\alpha}{\mathfrak{p}}\right)_m = 0$.

$\mathrm{(2)}$ If $\alpha\notin\mathfrak{p}$, then $\left(\frac{\alpha}{\mathfrak{p}}\right)_m =
\zeta_m^i$ with $i\in\mathbb{Z}$, where $\zeta_m^i$
is the unique $m$-th root of unity in $D$ such that
$$\alpha^{(N(\mathfrak{p})-1)/m}\equiv \zeta_m^i\qquad (\textup{mod} \;\mathfrak{p}),$$
where N$(\mathfrak{p})$ is the absolute norm of the ideal $\mathfrak{p}$.

$\mathrm{(3)}$ If $\mathfrak{a}\subset D$ is an arbitrary ideal prime to $m$ and $\mathfrak{a}=\prod\mathfrak{p}_i^{n_i}$ is its factorization
 as a product of prime ideals, then
 $$\left(\frac{\alpha}{\mathfrak{a}}\right)_m=\prod\left(\frac{\alpha}{\mathfrak{p}_i}\right)_m^{n_i}.$$
We set $\left(\frac{\alpha}{D}\right)_m$ = 1.

 $\mathrm{(4)}$ If $\beta\in D$ and $\beta$ is prime to $m$ define
 $\left(\frac{\alpha}{\beta}\right)_m=\left(\frac{\alpha}{\beta D}\right)_m$.

We will need the following proposition which can be found in
\cite[Chapter 14, Corollary 2, p. 218]{ir}.

\smallskip

\bpro\label{propo1} Suppose $A, B\subset\mathbb{Z}[\zeta_m]$ are ideals prime to
$m$ and $A=(\alpha)$ is principal with gcd$(N(A), N(B))=1$. Then
$$\left(\frac{N(B)}{\alpha}\right)_m=\left(\frac{\varepsilon(\alpha)}{B}\right)_m\left(\frac{\alpha}{N(B)}\right)_m$$
where $\varepsilon(\alpha)=\pm\zeta_m^i$ for some $i\in\mathbb{Z}$. \epro

\smallskip

From the above proposition, we can obtain a certain special
$2p$-th reciprocity law.

\smallskip
\bpro\label{reciprocitylaw} Let $M\equiv1\pmod{4p^2}$ be an integer with $p$ an odd
prime. Let $\pi\in \mathbb{Z}[\zeta_{2p}]$ be coprime with $2pM$. Suppose
$M>1$ is a prime, then we have
$$\left(\frac{M}{\pi}\right)_{2p}=\left(\frac{\pi}{M}\right)_{2p}.$$\epro
\bpf{} Let $\mathfrak{P}$ be a prime ideal of
$\mathbb{Z}[\zeta_{2p}]$ lying over $M$. Since $M\equiv1$ (mod
$2p$), we have N$(\mathfrak{P})=M$. By Proposition \ref{propo1},
$\left(\frac{N(\mathfrak{P})}{\pi}\right)_{2p}=\left(\frac{\varepsilon(\pi)}{\mathfrak{P}}\right)_{2p}\left(\frac{\pi}{N(\mathfrak{P})}\right)_{2p}$
which implies
$\left(\frac{M}{\pi}\right)_{2p}=\left(\frac{\varepsilon(\pi)}{\mathfrak{P}}\right)_{2p}\left(\frac{\pi}{M}\right)_{2p}$.
And
$\left(\frac{\varepsilon(\pi)}{\mathfrak{P}}\right)_{2p}\equiv\varepsilon(\pi)^{(M-1)/2p}\equiv(\pm\zeta_{2p}^i)^{(M-1)/2p}=1$
(mod $\mathfrak{P}$), as $2p\mid\frac{M-1}{2p}$. Then $\left(\frac{\varepsilon(\pi)}{\mathfrak{P}}\right)_{2p}=1$
and $\left(\frac{M}{\pi}\right)_{2p}=\left(\frac{\pi}{M}\right)_{2p}.$
\qed
\epf
\smallskip

\section{Primality test for $M=(2p)^{2^n}+1$}\label{SecMain}

We first define the polynomials $G_n(x)(n\geq0)$ by $G_0(x)=1$, $G_1(x)=x$,
and for $n\geq2$ define $G_n(x)$ recursively by the
formulas:

$$G_n(x)=\left\{\begin{aligned}
&G_{(n-1)/2}(x)G_{(n+1)/2}(x)-x,\quad \textup{if} \; n \;\textup{is odd},\\
&G_{n/2}(x)^2-2,\qquad\qquad\;\;\;\;\textup{if} \;n\;\textup{is even}.
\end{aligned}
\right.$$

Clearly all the coefficients of $x^i(i\geq0)$ of $G_n(x)$ can be computed in $O(\textup{log}\;n)$
steps. It is easy to see that $G_n(x)(n\geq0)$ is a monic polynomial of degree $n$, and that
$G_n(x+x^{-1})=x^n+x^{-n}$ for $n>0$.

Let $$f^{(j)}(x_1,\ldots,x_{(p-1)/2})=\sum\limits_{1\leq
i_1<\cdots<i_j\leq (p-1)/2}x_{i_1}\cdots x_{i_j}, 1\leq j\leq
(p-1)/2$$ be the $j$-th elementary symmetric
polynomial of $x_1, \ldots, x_{(p-1)/2}$. Let
$D=\mathbb{Z}[\zeta_{2p}]$ be the ring of integers of the cyclotomic field $L =
\mathbb{Q}(\zeta_{2p})$.
Let $G=$Gal$(\mathbb{Q}(\zeta_{2p})/\mathbb{Q})$ be the Galois group
of $\mathbb{Q}(\zeta_{2p})$ over $\mathbb{Q}$. For every integer $c$ with gcd$(c,2p)=1$
denote by $\sigma_c$ the element of $G$ that sends $\zeta_{2p}$ to
$\zeta_{2p}^c$. We know that Gal$(\mathbb{Q}(\zeta_{2p})/\mathbb{Q})=\{\sigma_{\pm(2i-1)}\;|\;1\leq i\leq(p-1)/2\}$.
 For $\tau$ in $\mathbb{Z}[G]$ and $\alpha$ in $L$ with $\alpha\neq0$ we
often denote by $\alpha^\tau$ to the action of the element $\tau$ of
$\mathbb{Z}[G]$ on the element $\alpha$ of $L$, that is,
$$\alpha^{\tau}:=\prod\limits_{\sigma\in G}\sigma(\alpha)^{k_{\sigma}}, \;\textup{if}
\;\tau=\sum\limits_{\sigma\in G}k_{\sigma}\sigma\;\textup{where}\; k_{\sigma}\in\mathbb{Z}.$$
If $\tau\in G$, we
will either write $\alpha^\tau$ or $\tau(\alpha)$. We also write $\sigma_1=1$ in $\mathbb{Z}[G]$.

Let $K=\mathbb{Q}(\zeta_{2p}+\zeta_{2p}^{-1})$ be the maximal real
subfield of $L$. We know that  Gal$(K/\mathbb{Q})=\{\sigma_{2i-1}|_{K}\;|\;1\leq i\leq(p-1)/2\}$.
Let $\pi\in D$ with $\pi\notin \mathbb{R}$. We
denote
$\alpha=(\pi/\bar{\pi})^{\gamma}$ where $$\gamma=\sum\limits_{i=1}^{(p-1)/2}(2i-1)\sigma_{(2i-1)^{-1}}\in\mathbb{Z}[G]$$
and a bar indicates the complex conjugation, and $(2i-1)^{-1}$ is the inverse
of $2i-1$ in the multiplicative group $(\mathbb{Z}/2p\mathbb{Z})^*$. Obviously, we have $\alpha\overline{\alpha}=1$. Next we define
$(p-1)/2$ many sequences: $\{S_k^{(j)}|_{k\geq0}\}, \; 1\leq j\leq
(p-1)/2$ by
$S_k^{(j)}=f^{(j)}(\alpha_1^{(k)},\ldots,\alpha_{(p-1)/2}^{(k)})$, where
$\alpha_i^{(k)}=\sigma_{2i-1}(\alpha^{(2p)^k}+\bar{\alpha}^{(2p)^k}),\;i=1,\ldots,(p-1)/2$.

We set
$F(x)=\sum\limits_{k=0}^{(p-1)/2}
G_k(x):=\sum\limits_{j=0}^{(p-1)/2}(-1)^ja_jx^{(p-1)/2-j}\in\mathbb{Z}[x]$,
clearly $a_0=1$. Since
$$F(\zeta_p+\zeta_p^{-1})=1+\sum_{k=1}^{(p-1)/2}G_k(\zeta_p+\zeta_p^{-1})=1+\sum_{k=1}^{(p-1)/2}(\zeta_p^k+\zeta_p^{-k})=0$$
and $F(x)$ is a monic polynomial of degree $(p-1)/2$, so $F(x)$ is the minimal polynomial of $\zeta_p+\zeta_p^{-1}$ over $\mathbb{Q}$.

Our primality test for numbers
$M=(2p)^{2^n}+1$ is described as follows.

\smallskip

\bthm\label{maintheorem}  Let $S_k^{(j)}$ and $ a_j$ be as before. Let $M=(2p)^{2^n}+1$ with $n\geq1$, $p$ be an odd
prime and $r=2^n$. Let $\pi\in\mathbb{Z}[\zeta_{2p}]$ be coprime with $2pM$
such that $\pi\notin\mathbb{R}$ and  $\left(\frac{M}{\pi}\right)_{2p}\neq\pm1$.
Suppose that if $x^{p-1}\equiv1\pmod{p^r}$ and $1<x<p^r$
then $x$ does not divide $M$. Then $M$ is prime if and only if one
of the following holds:

$\mathrm{(i)}$ \;$\left(\frac{M}{\pi}\right)_{2p}=\zeta_p^l$ for some $l\in\mathbb{Z}$ and
$l\nequiv0$ $\textup{(mod}$ $p\textup{)}$, and $S_{r-1}^{(j)}\equiv a_j$ \textup{(mod} $M$\textup{)} for each
$1\leq j\leq(p-1)/2$;

$\mathrm{(ii)}$ $\left(\frac{M}{\pi}\right)_{2p}=-\zeta_p^l$ for some $l\in\mathbb{Z}$ and
$l\nequiv0$ \textup{(mod} $p$\textup{)}, and $S_{r-1}^{(j)}\equiv(-1)^ja_j$ \textup{(mod} $M$\textup{)} for each
  $1\leq j\leq(p-1)/2$. \ethm
 \bpf{}\;
We first show the necessity for the primality of $M$.  Suppose
then $M$ is a prime, since $\pi$ is prime to $2pM$,  applying
Proposition \ref{reciprocitylaw} we get $\left(\frac{M}{\pi}\right)_{2p}=\left(\frac{\pi}{M}\right)_{2p}$.
From $M\equiv 1$ (mod $2p$), the ideal $MD$ factors in $D$ as a
product of $p-1$ distinct prime ideals. We write
$$MD=(\mathfrak{p}\bar{\mathfrak{p}})^{\sum\limits_{i=1}^{(p-1)/2}\sigma_{2i-1}},$$
thus

\[
\begin{array}{ll}
\left(\frac{M}{\pi}\right)_{2p}=\left(\frac{\pi}{M}\right)_{2p}
&=\prod\limits_{i=1}^{(p-1)/2}\left(\frac{\pi}{(\mathfrak{p}\bar{\mathfrak{p}})^{\sigma_{2i-1}}}\right)_{2p}\\\\
&=\prod\limits_{i=1}^{(p-1)/2}\left(\frac{(\frac{\pi}{\bar{\pi}})^{(2i-1)\sigma_{(2i-1)^{-1}}}}{\mathfrak{p}}\right)_{2p}
=\left(\frac{(\frac{\pi}{\bar{\pi}})^{\sum\limits_{k=1}^{(p-1)/2}(2i-1)\sigma_{(2i-1)^{-1}}}}{\mathfrak{p}}\right)_{2p}\\\\
&=\left(\frac{\alpha}{\mathfrak{p}}\right)_{2p}\equiv\alpha^{(M-1)/2p}\equiv\alpha^{(2p)^{r-1}}
\quad \textup{(mod} \;\mathfrak{p}\textup{)}
\end{array}
\]
Since $\mathfrak{p}$ is an arbitrary prime ideal lying over $M$, we
have

$$\left(\frac{M}{\pi}\right)_{2p}\equiv\alpha^{(2p)^{r-1}}\quad\textup{( mod} \; M\textup{)}.$$
It implies
$$\alpha_1^{(r-1)}=\alpha^{(2p)^{r-1}}+\bar{\alpha}^{(2p)^{r-1}}\equiv\left(\frac{M}{\pi}\right)_{2p}
+\left(\frac{M}{\pi}\right)_{2p}^{-1}\quad \textup{(mod} \; M\textup{)},$$
and for all $1\leq i\leq(p-1)/2$, we obtain
$$\alpha_i^{(r-1)}=\sigma_{2i-1}(\alpha^{(2p)^{r-1}}+\bar{\alpha}^{(2p)^{r-1}})
\equiv\left(\frac{M}{\pi}\right)_{2p}^{2i-1}+\left(\frac{M}{\pi}\right)_{2p}^{1-2i}\quad
\textup{(mod} \; M\textup{)}.$$ Hence for all $1\leq j\leq(p-1)/2$,
\[
\begin{array}{ll}
S_{r-1}^{(j)}
&=f^{(j)}(\alpha_1^{(r-1)},\ldots,\alpha_{(p-1)/2}^{(r-1)})\\\\
&\equiv
f^{(j)}\left(\left(\frac{M}{\pi}\right)_{2p}+\left(\frac{M}{\pi}\right)_{2p}^{-1},\ldots,\left(\frac{M}{\pi}\right)_{2p}^{p-2}+
\left(\frac{M}{\pi}\right)_{2p}^{2-p}\right)
\quad \textup{(mod} \; M\textup{)}.
\end{array}
\]
$\mathrm{(i)}$ \;If $\left(\frac{M}{\pi}\right)_{2p}=\zeta_p^l$ for some $l$ with
$l\nequiv0$ $\textup{(mod}$ $p\textup{)}$. Since the minimal polynomial of $\zeta_p^l+\zeta_p^{-l}$ is
$F(x)$ defined before, we obtain for all $j:\;1\leq j\leq(p-1)/2$,
$a_j=f^{(j)}(\zeta_p+\zeta_p^{-1},\zeta_p^3+\zeta_p^{-3},\ldots,\zeta_p^{p-2}+\zeta_p^{2-p})$
because of
$F(x)=\prod\limits_{i=1}^{(p-1)/2}[x-(\zeta_p^{2i-1}+\zeta_p^{1-2i})]$.
Hence from above
$$S_{r-1}^{(j)}\equiv f^{(j)}(\zeta_p+\zeta_p^{-1},\zeta_p^3+\zeta_p^{-3},\ldots,\zeta_p^{p-2}+\zeta_p^{2-p})
=a_j\quad \textup{(mod} \; M\textup{)},$$ for all $j=1,\ldots,(p-1)/2$.

$\mathrm{(ii)}$ \;If $\left(\frac{M}{\pi}\right)_{2p}=-\zeta_p^l$ for some $l$ with
$l\nequiv0$ $\textup{(mod}$ $p\textup{)}$. Then by the property of elementary symmetric
polynomial, we have
\[
\begin{array}{ll}
S_{r-1}^{(j)}
&\equiv f^{(j)}\left(\left(\frac{M}{\pi}\right)_{2p}+\left(\frac{M}{\pi}\right)_{2p}^{-1},
\ldots,\left(\frac{M}{\pi}\right)_{2p}^{p-2}+\left(\frac{M}{\pi}\right)_{2p}^{2-p}\right)\\\\
&=f^{(j)}(-\zeta_p-\zeta_p^{-1},-\zeta_p^3-\zeta_p^{-3},\ldots,-\zeta_p^{p-2}-\zeta_p^{2-p})\\\\
&=(-1)^jf^{(j)}(\zeta_p+\zeta_p^{-1},\zeta_p^3+\zeta_p^{-3},\ldots,\zeta_p^{p-2}+\zeta_p^{2-p})=(-1)^ja_j
\quad \textup{(mod} \; M\textup{)},
\end{array}
\]
for all $j=1,\ldots,(p-1)/2$. This completes the proof of necessity.

Now we turn to the proof of sufficiency. Let $q$ be an arbitrary
prime divisor of $M$.  Let $\mathfrak{q}$ be a prime ideal in the
ring of integers of $K$ lying over $q$, and $\mathfrak{Q}$ be a
prime ideal of $D$ lying over $\mathfrak{q}$. Let
$\beta=\alpha^{(2p)^{r-1}}+\bar{\alpha}^{(2p)^{r-1}}\in K$, then
$S_{r-1}^{(j)}=f^{(j)}(\beta,\sigma_3(\beta),\ldots,\sigma_{p-2}(\beta))$,
$j=1,\ldots,(p-1)/2$.

$\mathrm{(i)}$ If $S_{r-1}^{(j)}\equiv a_j$ $\textup{(mod}$ $M\textup{)}$, then
$f^{(j)}(\beta,\sigma_3(\beta),\ldots,\sigma_{p-2}(\beta))\equiv
a_j$ $\textup{(mod}$ $\mathfrak{q}\textup{)}$, hence we have
\[
\begin{array}{ll}
0&=(\beta-\beta)(\beta-\sigma_3(\beta))\ldots(\beta-\sigma_{p-2}(\beta))\\\
&=\beta^{(p-1)/2}+\sum\limits_{j=1}^{(p-1)/2}(-1)^jf^{(j)}(\beta,\sigma_3(\beta),\ldots,\sigma_{p-2}(\beta))
\beta^{(p-1)/2-j}\\\\
&\equiv\beta^{(p-1)/2}+\sum\limits_{j=1}^{(p-1)/2}(-1)^j
a_j\beta^{(p-1)/2-j}=F(\beta)\quad \textup{(mod} \; \mathfrak{q}\textup{)},
\end{array}
\]
and
\[
\begin{array}{ll}
0&\equiv
F(\alpha^{(2p)^{r-1}}+\bar{\alpha}^{(2p)^{r-1}})\\\\
&=1+\sum\limits_{k=1}^{(p-1)/2}\left[(\alpha^{(2p)^{r-1}})^k+(\bar{\alpha}^{(2p)^{r-1}})^k\right]
\quad \textup{(mod} \; \mathfrak{Q}\textup{)}.
\end{array}
\]
Multiplying both sides of the above congruence by
$\alpha^{(2p)^{r-1}\cdot(p-1)/2}=\bar{\alpha}^{-(2p)^{r-1}\cdot(p-1)/2}$
gives
$$\sum\limits_{k=0}^{p-1}\alpha^{(2p)^{r-1}k}\equiv0\quad \textup{(mod} \; \mathfrak{Q}\textup{)}.$$
It implies that the image of $\alpha^{(2p)^{r-1}}$ has order $p$ in
the multiplicative group $(D/\mathfrak{Q})^*$, and so the image of
$\alpha^{2^{r-1}}$ has order $p^r$.  This multiplicative group has
order N$(\mathfrak{Q})-1$ which divides $q^{p-1}-1$, $i.e.$,
$q^{p-1}\equiv1$ $\textup{(mod}$ $p^r\textup{)}$.  By the assumption $M$ is not
divisible by all solutions of equation $x^{p-1}\equiv1$ $\textup{(mod}$ $p^r\textup{)}$
 between $1$ and $p^r$.  Then $q>p^r>\sqrt{(2p)^r+1}=\sqrt{M}$,
 thus $q >\sqrt{M}$ for arbitrary prime
 divisor $q$ of $M$, that is to say $M$ is prime.

$\mathrm{(ii)}$ If $S_{r-1}^{(j)}\equiv(-1)^j a_j$ $\textup{(mod}$ $M\textup{)}$, then
$f^{(j)}(\beta,\sigma_3(\beta),\ldots,\sigma_{p-2}(\beta))\equiv
(-1)^j a_j$ $\textup{(mod}$ $\mathfrak{q}\textup{)}$, hence we get
\[
\begin{array}{ll}
0&=\beta^{(p-1)/2}+\sum\limits_{j=1}^{(p-1)/2}(-1)^jf^{(j)}(\beta,\sigma_3(\beta),\ldots,\sigma_{p-2}(\beta))
\beta^{(p-1)/2-j}\\\\
&\equiv\beta^{(p-1)/2}+\sum\limits_{j=1}^{(p-1)/2}
a_j\beta^{(p-1)/2-j}=(-1)^{(p-1)/2}F(-\beta)\quad \textup{(mod} \; \mathfrak{q}\textup{)},
\end{array}
\]
and
\[
\begin{array}{ll}
0&\equiv F(-\alpha^{(2p)^{r-1}}-\bar{\alpha}^{(2p)^{r-1}})\\\\
&=1+\sum\limits_{k=1}^{(p-1)/2}\left[(-\alpha^{(2p)^{r-1}})^k+(-\bar{\alpha}^{(2p)^{r-1}})^k\right]
\quad \textup{(mod} \; \mathfrak{Q}\textup{)}.
\end{array}
\]
Again multiplying both sides of the above congruence by
$\alpha^{(2p)^{r-1}\cdot(p-1)/2}=\bar{\alpha}^{-(2p)^{r-1}\cdot(p-1)/2}$
gives
$$\sum\limits_{k=0}^{p-1}(-1)^{k-(p-1)/2}\alpha^{(2p)^{r-1}k}\equiv0\quad
\textup{(mod} \; \mathfrak{Q}\textup{)}.$$ That is,
$$\sum\limits_{k=0}^{p-1}(-1)^k(\alpha^{(2p)^{r-1}})^k\equiv0\quad
\textup{(mod} \; \mathfrak{Q}\textup{)}.$$ Hence we obtain the image of
$\alpha^{(2p)^{r-1}}$ has order $2p$ in the multiplicative group
$(D/\mathfrak{Q})^*$, and so the image of $\alpha$ has order
$(2p)^r$. The same as above we have $(2p)^r$ divides the order of
group $(D/\mathfrak{Q})^*$, and $q^{p-1}\equiv1$ $\textup{(mod}$ $p^r\textup{)}$. By the
assumption we similarly get $q>p^r>\sqrt{(2p)^r+1}=\sqrt{M}$, and $q
>\sqrt{M}$ for arbitrary prime divisor $q$ of $M$, then $M$ is prime.
This completes the proof of sufficiency.
\qed\epf
%\smallskip

\textbf{Remark.}\quad $\mathrm{(i)}$ \;The methods on how to find $\pi$ and how to solve the equation $x^{p-1}\equiv1\pmod{p^r}$ involved in Theorem \ref{maintheorem} can be found in \cite{dl}. In general, this is not a difficult matter.

$\mathrm{(ii)}$ \; The testing sequences $\{S_k^{(j)}|{k\geq0}\}$, $1\leq j\leq(p-1)/2$ are sequences of
rational numbers. Hence the  computation of the sufficient and necessary condition in Theorem \ref{maintheorem}
is taken in residue class ring $\mathbb{Z}/M\mathbb{Z}$. And this benefits the explicit realization of
the primality test for $M$ in the algorithm.

\section{Primality tests for $p\leq19$}\label{SecExpli}

We know from \cite[Chapter 11]{wa} that $\mathbb{Z}[\zeta_{2p}]$ is a PID for $p\leq19$.
In this section we apply Theorem \ref{maintheorem} to the cases $3\leq p\leq19$ with
$p$ prime. By the definition of $G_k(x)$ in
 the previous section, we can compute $G_k(x)$, $0\leq k\leq9$, as follows:

$G_0(x)=1$, $G_1(x)=x$, $G_2(x)=x^2-2$, $G_3(x)=x^3-3x$, $G_4(x)=x^4-4x^2+2$,
$G_5(x)=x^5-5x^3+5x$, $G_6(x)=x^6-6x^4+9x^2-2$, $G_7(x)=x^7-7x^5+14x^3-7x$,
$G_8(x)=x^8-8x^6+20x^4-16x^2+2$, $G_9(x)=x^9-9x^7+27x^5-30x^3+9x$.

And we denote $F_1(x)=G_0(x)+G_1(x)=x+1$,
 $F_2(x)=G_0(x)+G_1(x)+G_2(x)=x^2+x-1$,
  $F_3(x)=\sum\limits_{k=0}^3G_k(x)=x^3+x^2-2x-1$,
  $F_4(x)=\sum\limits_{k=0}^5G_k(x)=x^5+x^4-4x^3-3x^2+3x+1$,
   $F_5(x)=\sum\limits_{k=0}^6G_k(x)=x^6+x^5-5x^4-4x^3+6x^2+3x-1$,
    $F_6(x)=\sum\limits_{k=0}^8G_k(x)=x^8+x^7-7x^6-6x^5+15x^4+10x^3-10x^2-4x+1$,
$F_7(x)=\sum\limits_{k=0}^9G_k(x)=x^9+x^8-8x^7-7x^6+21x^5+15x^4-20x^3-10x^2+5x+1$.
Next we can obtain $7$ special primality tests.
%\smallskip

\bpro\label{1} Let $M=6^{2^n}+1$, $n\geq1$ and $r=2^n$. Let
$\pi=2+3\zeta_3\in\mathbb{Z}[\zeta_6]$, and
$\alpha=\pi/\bar{\pi}$. We define sequence
$\{S_k\}$ with $S_0=\alpha+\bar{\alpha}$ and
$S_{k+1}=S_k^6-6S_k^4+9S_k^2-2$ for $k\geq0$. Then $M$ is prime if
and only if $S_{r-1}\equiv-1$ $\textup{(mod}$ $M\textup{)}$. \epro
\bpf{}\;
We let $L=\mathbb{Q}(\zeta_6)$, then
$Norm_{L/\mathbb{Q}}(\pi)=\pi\bar{\pi}=(2+3\zeta_3)(-1-3\zeta_3)=7$.
Since $M\equiv2$ $\textup{(mod}$ $7\textup{)}$, then we have $\left(\frac{M}{\pi}\right)_6\equiv
M^{(7-1)/6}=M\equiv2\equiv\zeta_3^2$ $\textup{(mod}$ $\pi\textup{)}$,
and $\left(\frac{M}{\pi}\right)_6=\zeta_3^2$. Let
$S_k=\alpha^{6^k}+\bar{\alpha}^{6^k}$, $k\geq0$. We can easily verify that
$S_k$ satisfies the recurrent relation in the assumption. We use the same polynomial
$F(x)$ as in the above section. Here
$F(x)=F_1(x)=x+1$ implies $a_1=-1$. Hence by the necessity of
Theorem \ref{maintheorem} we obtain that if $M$ is prime then $S_{r-1}\equiv-1$
$\textup{(mod}$ $M\textup{)}$. This completes the proof of necessity.

And by the proof of sufficiency of Theorem \ref{maintheorem}, we get
if $S_{r-1}\equiv-1$ $\textup{(mod}$ $M\textup{)}$, then $3^r$ divides $q^2-1$ for
arbitrary prime divisor $q$ of $M$, $i.e.$, $3^r$ divides one of $q+1$ and
$q-1$ because of gcd$(q+1,q-1)=2$. So
$q\geq3^r-1>\sqrt{6^{r}+1}=\sqrt{M}$ and we get $M$ is prime. This
completes the proof of sufficiency.
\qed\epf
%\smallskip

\textbf{Remark.}\quad This primality test is explicit. Compare to the primality test of $G_n=6^{2^n}+1$ in \cite{wi} they own
considerable complexity of running time which is polynomial in log$_{2}(M)$. And it seems that our test
is more succinct.

\bpro\label{2} Let $M=10^{2^n}+1$, $n\geq1$ and $r=2^n$. Let
$\pi=1-\zeta_5-\zeta_5^3\in\mathbb{Z}[\zeta_{10}]$,
$\alpha=(\pi/\bar{\pi})^{1+3\sigma_{-3}}$. We define sequences $\{S_k^{(1)}\}$ and $\{S_k^{(2)}\}$ with
$S_k^{(1)}=\alpha_1^{(k)}+\alpha_2^{(k)}$, $S_k^{(2)}=\alpha_1^{(k)}\cdot\alpha_2^{(k)}$,
$k\geq0$, where
$\alpha_1^{(k)}=\alpha^{10^k}+\bar{\alpha}^{10^k}$,
$\alpha_2^{(k)}=\sigma_3(\alpha_1^{(k)})$. Suppose that if $x^4\equiv1$
$\textup{(mod}$ $5^r\textup{)}$ and $1<x<5^r$ then $x$ does not divide $M$. Then $M$ is
prime if and only if $S_{r-1}^{(1)}\equiv1\equiv-S_{r-1}^{(2)}$ $\textup{(mod}$ $M\textup{)}$.\epro
\bpf{}\;
We let $L=\mathbb{Q}(\zeta_{10})$ , then
 $Norm_{L/\mathbb{Q}}(\pi)=(\pi\bar{\pi})^{1+\sigma_3}=11$. Since
$M\equiv2$ $\textup{(mod}$ $11\textup{)}$, $n\geq1$, then
$\left(\frac{M}{\pi}\right)_{10}\equiv
M^{(11-1)/10}=M\equiv2\equiv-\zeta_5$ $\textup{(mod}$ $\pi\textup{)}$, and
$\left(\frac{M}{\pi}\right)_{10}=-\zeta_5$. We notice that here
$F(x)=F_2(x)=x^2+x-1$, which implies $a_1=-1$, $a_2=-1$.
Hence by the necessity of Theorem \ref{maintheorem} we obtain that if $M$ is
prime then $S_{r-1}^{(1)}\equiv-a_1=1$ $\textup{(mod}$ $M\textup{)}$,
$S_{r-1}^{(2)}\equiv a_2=-1$ $\textup{(mod}$ $M\textup{)}$. This completes the proof of necessity.

 By the sufficiency of Theorem \ref{maintheorem}, we easily get if
$S_{r-1}^{(1)}\equiv1\equiv-S_{r-1}^{(2)}$ $\textup{(mod}$ $M\textup{)}$,
i.e.,  $S_{r-1}^{(1)}\equiv-a_1$ $\textup{(mod}$ $M\textup{)}$,
$S_{r-1}^{(2)}\equiv a_2$ $\textup{(mod}$ $M\textup{)}$, then by assumption we have $M$ is
prime. This
completes the proof of sufficiency.
\qed\epf
%\smallskip

\bpro Let $M=14^{2^n}+1$, $n>1$ and $r=2^n$. Let
$\pi=1-\zeta_7+\zeta_7^4\in\mathbb{Z}[\zeta_{14}]$,
$\alpha=(\pi/\bar{\pi})^{1+3\sigma_5+5\sigma_3}$. We define
the sequences $\{S_k^{(1)}\}$, $\{S_k^{(2)}\}$ and $\{S_k^{(3)}\}$ with
$S_k^{(1)}=\alpha_1^{(k)}+\alpha_2^{(k)}+\alpha_3^{(k)}$,
$S_k^{(2)}=f^{(2)}(\alpha_1^{(k)},\alpha_2^{(k)},\alpha_3^{(k)})$,
$S_k^{(3)}=\alpha_1^{(k)}\alpha_2^{(k)}\alpha_3^{(k)}$, $k\geq0$, where
$\alpha_1^{(k)}=\alpha^{14^k}+\bar{\alpha}^{14^k}$,
$\alpha_2^{(k)}=\sigma_3(\alpha_1^{(k)})$,
$\alpha_3^{(k)}=\sigma_5(\alpha_1^{(k)})$. Suppose that if $x^6\equiv1$
$\textup{(mod}$ $7^r\textup{)}$ and $1<x<7^r$ then $x$ does not divide $M$. Then $M$ is
prime if and only if one of the following
holds:

$\mathrm{(i)}$ $M\equiv\pm8$ $\textup{(mod}$ $29\textup{)}$ and
$S_{r-1}^{(1)}\equiv1\equiv-S_{r-1}^{(3)}$ $\textup{(mod}$ $M\textup{)}$, $S_{r-1}^{(2)}\equiv-2$
$\textup{(mod}$ $M\textup{)}$;

$\mathrm{(ii)}$ $M\equiv-5$ $\textup{(mod}$ $29\textup{)}$ and
$S_{r-1}^{(1)}\equiv-1\equiv-S_{r-1}^{(3)}$ $\textup{(mod}$ $M\textup{)}$, $S_{r-1}^{(2)}\equiv-2$
$\textup{(mod}$ $M\textup{)}$.\epro

\bpf{}\; We let $L=\mathbb{Q}(\zeta_{14})$ , then
$Norm_{L/\mathbb{Q}}(\pi)=(\pi\bar{\pi})^{1+\sigma_3+\sigma_5}=29$. Since
$M\equiv\pm8,-5$ $\textup{(mod}$ $29\textup{)}$ , $n>1$, we have
$\left(\frac{M}{\pi}\right)_{14}\equiv
M^{(29-1)/14}=M^2\equiv6,-4\equiv-\zeta_7^3, \;\zeta_7$ $\textup{(mod}$ $\pi\textup{)}$,
 and $\left(\frac{M}{\pi}\right)_{14}=-\zeta_7^3, \; \zeta_7$. We notice
that here $F(x)=F_3(x)=x^3+x^2-2x-1$, which implies
$a_1=-1$, $a_2=-2$, $a_3=1$. Hence by the necessity of  Theorem
\ref{maintheorem}, if $M$ is prime and $M\equiv\pm8$ $\textup{(mod}$ $29\textup{)}$, then we obtain
$S_{r-1}^{(1)}\equiv-a_1=1$ $\textup{(mod}$ $M\textup{)}$, $S_{r-1}^{(2)}\equiv a_2=-2$
$\textup{(mod}$ $M\textup{)}$, $S_{r-1}^{(3)}\equiv-a_3=-1$ $\textup{(mod}$ $M\textup{)}$
because of $\left(\frac{M}{\pi}\right)_{14}=-\zeta_7^3$. If  $M$ is
prime and $M\equiv-5$ $\textup{(mod}$ $29\textup{)}$, then we similarly get
$S_{r-1}^{(1)}\equiv-1\equiv-S_{r-1}^{(3)}$ $\textup{(mod}$ $M\textup{)}$, $S_{r-1}^{(2)}\equiv-2$
$\textup{(mod}$ $M\textup{)}$ for $\left(\frac{M}{\pi}\right)_{14}=\zeta_7$. This completes the proof of necessity.

By the sufficiency of  Theorem \ref{maintheorem}, whatever $\mathrm{(i)}$ or
$\mathrm{(ii)}$ holds, we can always deduce that $M$ is a prime. This
completes the proof of sufficiency.
\qed\epf
%\smallskip

\bpro Let $M=22^{2^n}+1$, $n\geq1$ and $r=2^n$. Let
$\pi=1+\zeta_{11}^7+\zeta_{11}^8\in\mathbb{Z}[\zeta_{22}]$,
$\alpha=(\pi/\bar{\pi})^{\tau}$, where $\tau=1+3\sigma_{-7}+5\sigma_9+7\sigma_{-3}+9\sigma_5$.
We define the sequences $\{S_k^{(j)}\}$, $1\leq j\leq5$ with
$S_k^{(j)}=f^{(j)}(\alpha_1^{(k)},\ldots,\alpha_5^{(k)})$, $k\geq0$,
where $\alpha_1^{(k)}=\alpha^{22^k}+\bar{\alpha}^{22^k}$, and
$\alpha_i^{(k)}=\sigma_{2i-1}(\alpha_1^{(k)})$, $2\leq i\leq5$. Suppose
that if $x^{10}\equiv1$ $\textup{(mod}$ $11^r\textup{)}$ and $1<x<11^r$ then $x$ does
not divide $M$. Then $M$ is prime if and only if
$S_{r-1}^{(1)}\equiv-1\equiv S_{r-1}^{(5)}$ $\textup{(mod}$ $M\textup{)}$, $S_{r-1}^{(2)}\equiv-4$
$\textup{(mod}$ $M\textup{)}$ and $S_{r-1}^{(3)}\equiv3\equiv S_{r-1}^{(4)}$ $\textup{(mod}$ $M\textup{)}$ .
 \epro

 \bpf{}\;We let $L=\mathbb{Q}(\zeta_{22})$ , then
$Norm_{L/\mathbb{Q}}(\pi)=(\pi\bar{\pi})^{\sum\limits_{i=1}^5\sigma_{2i-1}}=23$. Since
$M\equiv2$ $\textup{(mod}$ $23\textup{)}$, $n\geq1$, we have
$\left(\frac{M}{\pi}\right)_{22}\equiv
M^{(23-1)/22}=M\equiv2\equiv\zeta_{11}^2$ $\textup{(mod}$ $\pi\textup{)}$, and
$\left(\frac{M}{\pi}\right)_{22}=\zeta_{11}^2$. Also we notice that here
$F(x)=F_4(x)=x^5+x^4-4x^3-3x^2+3x+1$, which implies $a_1=-1$,
$a_2=-4$, $a_3=3$, $a_4=3$, $a_5=-1$. Hence by the necessity of Theorem \ref{maintheorem}, we have if
$M$ is prime, then
$S_{r-1}^{(1)}\equiv a_1=-1$ $\textup{(mod}$ $M\textup{)}$,
$S_{r-1}^{(2)}\equiv a_2=-4$ $\textup{(mod}$ $M\textup{)}$,
$S_{r-1}^{(3)}\equiv a_3=3$ $\textup{(mod}$ $M\textup{)}$,
$S_{r-1}^{(4)}\equiv a_4=3$ $\textup{(mod}$ $M\textup{)}$,
$S_{r-1}^{(5)}\equiv a_5=-1$ $\textup{(mod}$ $M\textup{)}$.
 This completes the proof of necessity.

 By the sufficiency
of Theorem \ref{maintheorem} and the assumption, clearly we have $M$ is a prime. This completes
 the proof of sufficiency.
\qed\epf
%\smallskip

\bpro Let $M=26^{2^n}+1$, $n>1$ and $r=2^n$. Let
$\pi=1+\zeta_{13}^2+\zeta_{13}^5\in\mathbb{Z}[\zeta_{26}]$,
$\alpha=(\pi/\bar{\pi})^{\tau}$, where $\tau=1+3\sigma_9+5\sigma_{-5}+7\sigma_{-11}+9\sigma_3+11\sigma_{-7}$.
We define the sequences $\{S_k^{(j)}\}$, $1\leq j\leq6$ with
$S_k^{(j)}=f^{(j)}(\alpha_1^{(k)},\ldots,\alpha_6^{(k)})$, $k\geq0$,
where $\alpha_1^{(k)}=\alpha^{26^k}+\bar{\alpha}^{26^k}$, and
$\alpha_i^{(k)}=\sigma_{2i-1}(\alpha_1^{(k)})$, $2\leq i\leq6$. Suppose
that if $x^{12}\equiv1$ $\textup{(mod}$ $13^r\textup{)}$ and $1<x<13^r$ then $x$ does
not divide $M$. Then $M$ is prime if and only if one of the following holds:

$\mathrm{(i)}$ $M\equiv25,\pm16,-6,11,-24,-5,-10,17$ $\textup{(mod}$ $53\textup{)}$ and
$S_{r-1}^{(1)}\equiv-1\equiv S_{r-1}^{(6)}$ $\textup{(mod}$ $M\textup{)}$,
$S_{r-1}^{(2)}\equiv-5$ $\textup{(mod}$ $M\textup{)}$,
$S_{r-1}^{(3)}\equiv4$ $\textup{(mod}$ $M\textup{)}$,
$S_{r-1}^{(4)}\equiv6$ $\textup{(mod}$ $M\textup{)}$ and
$S_{r-1}^{(5)}\equiv-3$ $\textup{(mod}$ $M\textup{)}$;

$\mathrm{(ii)}$ $M\equiv14,-8,-3$ $\textup{(mod}$ $53\textup{)}$ and
$S_{r-1}^{(1)}\equiv1\equiv-S_{r-1}^{(6)}$ $\textup{(mod}$ $M\textup{)}$,
$S_{r-1}^{(2)}\equiv-5$ $\textup{(mod}$ $M\textup{)}$,
$S_{r-1}^{(3)}\equiv-4$ $\textup{(mod}$ $M\textup{)}$,
$S_{r-1}^{(4)}\equiv6$ $\textup{(mod}$ $M\textup{)}$ and
$S_{r-1}^{(5)}\equiv3$ $\textup{(mod}$ $M\textup{)}$.\epro

\bpf{}\;We let $L=\mathbb{Q}(\zeta_{26})$ , then
$Norm_{L/\mathbb{Q}}(\pi)=(\pi\bar{\pi})^{\sum\limits_{i=1}^6\sigma_{2i-1}}=53$. Since
$M\equiv25,\pm16,-6,\\
11,-24,-5,-10,17,14,-8,-3$ $\textup{(mod}$ $53\textup{)}$, $n>1$, we have
$\left(\frac{M}{\pi}\right)_{26}\equiv
M^{(53-1)/26}=M^2\equiv-11,-9,-17,15,-7,25,-6,24,-16,11,9\equiv\zeta_{13}^3,\;
\zeta_{13}^4,\;\zeta_{13}^5,\;\zeta_{13}^6,\;\zeta_{13}^7,\;\zeta_{13}^8,\;
\zeta_{13}^9,\;\zeta_{13}^{10},\;-\zeta_{13}^2,\;-\zeta_{13}^3,\\
-\zeta_{13}^4$ $\textup{(mod}$ $\pi\textup{)}$, and
$\left(\frac{M}{\pi}\right)_{26}=\zeta_{13}^3,\;
\zeta_{13}^4,\;\zeta_{13}^5,\;\zeta_{13}^6,\;\zeta_{13}^7,\;\zeta_{13}^8,\;
\zeta_{13}^9,\;\zeta_{13}^{10},\;-\zeta_{13}^2,\;-\zeta_{13}^3,
\;-\zeta_{13}^4$. Here we notice that
$F(x)=F_5(x)=x^6+x^5-5x^4-4x^3+6x^2+3x-1$, which implies $a_1=-1$,
$a_2=-5$, $a_3=4$, $a_4=6$, $a_5=-3$, $a_6=-1$. Hence by the necessity of Theorem \ref{maintheorem}, if
$M$ is prime and $M\equiv25,\pm16,-6,11,-24,-5,-10,17$ $\textup{(mod}$ $53\textup{)}$ then
we obtain
$S_{r-1}^{(1)}\equiv a_1=-1$ $\textup{(mod}$ $M\textup{)}$,
$S_{r-1}^{(2)}\equiv a_2=-5$ $\textup{(mod}$ $M\textup{)}$,
$S_{r-1}^{(3)}\equiv a_3=4$ $\textup{(mod}$ $M\textup{)}$,
$S_{r-1}^{(4)}\equiv a_4=6$ $\textup{(mod}$ $M\textup{)}$,
$S_{r-1}^{(5)}\equiv a_5=-3$ $\textup{(mod}$ $M\textup{)}$,
$S_{r-1}^{(6)}\equiv a_6=-1$ $\textup{(mod}$ $M\textup{)}$,
due to here $\left(\frac{M}{\pi}\right)_{26}$ is a primitive $13$-th root of unity.
If $M$ is prime and  $M\equiv14,-8,-3$ $\textup{(mod}$ $53\textup{)}$, then we also can get
$S_{r-1}^{(1)}\equiv1\equiv-S_{r-1}^{(6)}$ $\textup{(mod}$ $M\textup{)}$,
$S_{r-1}^{(2)}\equiv-5$ $\textup{(mod}$ $M\textup{)}$,
$S_{r-1}^{(3)}\equiv-4$ $\textup{(mod}$ $M\textup{)}$,
$S_{r-1}^{(4)}\equiv6$ $\textup{(mod}$ $M\textup{)}$,
$S_{r-1}^{(5)}\equiv3$ $\textup{(mod}$ $M\textup{)}$
because $\left(\frac{M}{\pi}\right)_{26}$ is a primitive $26$-th root of unity.
This completes the proof of necessity.

 Next by the sufficiency
of Theorem \ref{maintheorem} and the assumption, whatever $\mathrm{(i)}$ or
$\mathrm{(ii)}$ holds, we can easily deduce that $M$ is a prime. This
completes the proof of sufficiency.
\qed\epf
%\smallskip

\bpro Let $M=34^{2^n}+1$, $n\geq1$ and $r=2^n$. Let
$\pi=1+\zeta_{17}^2+\zeta_{17}^9\in\mathbb{Z}[\zeta_{34}]$,
$\alpha=(\pi/\bar{\pi})^{\tau}$, where $\tau=1+3\sigma_{-11}+5\sigma_7+7\sigma_5+9\sigma_{-15}+11\sigma_{-3}+
13\sigma_{-13}+15\sigma_{-9}$.
We define the sequences $\{S_k^{(j)}\}$, $1\leq j\leq8$ with
$S_k^{(j)}=f^{(j)}(\alpha_1^{(k)},\ldots,\alpha_8^{(k)})$, $k\geq0$,
where $\alpha_1^{(k)}=\alpha^{34^k}+\bar{\alpha}^{34^k}$, and
$\alpha_i^{(k)}=\sigma_{2i-1}(\alpha_1^{(k)})$, $2\leq i\leq8$. Suppose
that if $x^{16}\equiv1$ $\textup{(mod}$ $17^r\textup{)}$ and $1<x<17^r$ then $x$ does
not divide $M$. Then $M$ is prime if and only if one of the following holds:

$\mathrm{(i)}$ $M\equiv-21,15$ $\textup{(mod}$ $103\textup{)}$ and
$S_{r-1}^{(1)}\equiv-1\equiv-S_{r-1}^{(8)}$ $\textup{(mod}$ $M\textup{)}$,
$S_{r-1}^{(2)}\equiv-7$ $\textup{(mod}$ $M\textup{)}$,
$S_{r-1}^{(3)}\equiv6$ $\textup{(mod}$ $M\textup{)}$,
$S_{r-1}^{(4)}\equiv15$ $\textup{(mod}$ $M\textup{)}$,
$S_{r-1}^{(5)}\equiv-10\equiv S_{r-1}^{(6)}$ $\textup{(mod}$ $M\textup{)}$ and
$S_{r-1}^{(7)}\equiv4$ $\textup{(mod}$ $M\textup{)}$;

$\mathrm{(ii)}$ $M\equiv35,24,-2,-9,10,-30$ $\textup{(mod}$ $103\textup{)}$ and
$S_{r-1}^{(1)}\equiv1\equiv S_{r-1}^{(8)}$ $\textup{(mod}$ $M\textup{)}$,
$S_{r-1}^{(2)}\equiv-7$ $\textup{(mod}$ $M\textup{)}$,
$S_{r-1}^{(3)}\equiv-6$ $\textup{(mod}$ $M\textup{)}$,
$S_{r-1}^{(4)}\equiv15$ $\textup{(mod}$ $M\textup{)}$,
$S_{r-1}^{(5)}\equiv10\equiv-S_{r-1}^{(6)}$ $\textup{(mod}$ $M\textup{)}$ and
$S_{r-1}^{(7)}\equiv-4$ $\textup{(mod}$ $M\textup{)}$.\epro

\bpf{}\;We let $L=\mathbb{Q}(\zeta_{34})$ , then
$Norm_{L/\mathbb{Q}}(\pi)=(\pi\bar{\pi})^{\sum\limits_{i=1}^8\sigma_{2i-1}}=103$. Since
$M\equiv-21,15,35,24,-2,-9,10,-30$ $\textup{(mod}$ $103\textup{)}$, $n\geq1$, we have
$\left(\frac{M}{\pi}\right)_{34}\equiv
M^{(103-1)/34}=M^3\equiv9,-24,27,22,-8,-30,-14\equiv\zeta_{17}^2,\;
\zeta_{17}^7,\;-\zeta_{17}^3,\;-\zeta_{17}^4,\;-\zeta_{17}^6,\;-\zeta_{17}^{10},\;
-\zeta_{17}^{13}$ $\textup{(mod}$ $\pi\textup{)}$, where $(-2)^3\equiv(-9)^3\equiv-8$
$\textup{(mod}$ $103\textup{)}$, thus
$\left(\frac{M}{\pi}\right)_{34}=\zeta_{17}^2,\;
\zeta_{17}^7,\;-\zeta_{17}^3,\;-\zeta_{17}^4,\;-\zeta_{17}^6,\;-\zeta_{17}^{10},\;
-\zeta_{17}^{13}$. Again we notice that
$F(x)=F_6(x)=x^8+x^7-7x^6-6x^5+15x^4+10x^3-10x^2-4x+1$, which implies $a_1=-1$,
$a_2=-7$, $a_3=6$, $a_4=15$, $a_5=-10$, $a_6=-10$, $a_7=4$, $a_8=1$. Hence by the necessity of Theorem \ref{maintheorem}, if
$M$ is prime and $M\equiv-21,15$ $\textup{(mod}$ $103\textup{)}$ then
we obtain
$S_{r-1}^{(1)}\equiv a_1=-1$ $\textup{(mod}$ $M\textup{)}$,
$S_{r-1}^{(2)}\equiv a_2=-7$ $\textup{(mod}$ $M\textup{)}$,
$S_{r-1}^{(3)}\equiv a_3=6$ $\textup{(mod}$ $M\textup{)}$,
$S_{r-1}^{(4)}\equiv a_4=15$ $\textup{(mod}$ $M\textup{)}$,
$S_{r-1}^{(5)}\equiv a_5=-10$ $\textup{(mod}$ $M\textup{)}$,
$S_{r-1}^{(6)}\equiv a_6=-10$ $\textup{(mod}$ $M\textup{)}$,
$S_{r-1}^{(7)}\equiv a_7=4$ $\textup{(mod}$ $M\textup{)}$,
$S_{r-1}^{(8)}\equiv a_8=1$ $\textup{(mod}$ $M\textup{)}$,
since here $\left(\frac{M}{\pi}\right)_{34}$ is a primitive $17$-th root of unity.
If $M$ is prime and  $M\equiv35,24,-2,-9,10,-30$ $\textup{(mod}$ $103\textup{)}$, then we also obtain
$S_{r-1}^{(1)}\equiv1\equiv S_{r-1}^{(8)}$ $\textup{(mod}$ $M\textup{)}$,
$S_{r-1}^{(2)}\equiv-7$ $\textup{(mod}$ $M\textup{)}$,
$S_{r-1}^{(3)}\equiv-6$ $\textup{(mod}$ $M\textup{)}$,
$S_{r-1}^{(4)}\equiv15$ $\textup{(mod}$ $M\textup{)}$,
$S_{r-1}^{(5)}\equiv10\equiv-S_{r-1}^{(6)}$ $\textup{(mod}$ $M\textup{)}$,
$S_{r-1}^{(7)}\equiv-4$ $\textup{(mod}$ $M\textup{)}$,
because $\left(\frac{M}{\pi}\right)_{34}$ is a primitive $34$-th root of unity.
This completes the proof of necessity.

Similarly by the sufficiency
of Theorem \ref{maintheorem} and the assumption, whatever $\mathrm{(i)}$ or
$\mathrm{(ii)}$ holds, we always deduce that $M$ is a prime. This
completes the proof of sufficiency.
\qed\epf
%\smallskip

\bpro Let $M=38^{2^n}+1$, $n>1$ and $r=2^n$. Let
$\pi=-1-\zeta_{19}^2+\zeta_{19}^{15}\in\mathbb{Z}[\zeta_{38}]$,
$\alpha=(\pi/\bar{\pi})^{\tau}$, where $\tau=1+3\sigma_{13}+5\sigma_{-15}+7\sigma_{11}+9\sigma_{17}+11\sigma_7+
13\sigma_3+15\sigma_{-5}+17\sigma_9$.
We define the sequences $\{S_k^{(j)}\}$, $1\leq j\leq9$ with
$S_k^{(j)}=f^{(j)}(\alpha_1^{(k)},\ldots,\alpha_9^{(k)})$, $k\geq0$,
where $\alpha_1^{(k)}=\alpha^{38^k}+\bar{\alpha}^{38^k}$, and
$\alpha_i^{(k)}=\sigma_{2i-1}(\alpha_1^{(k)})$, $2\leq i\leq9$. Suppose
that if $x^{18}\equiv1$ $\textup{(mod}$ $19^r\textup{)}$ and $1<x<19^r$ then $x$ does
not divide $M$. Then $M$ is prime if and only if one of the following holds:

$\mathrm{(i)}$ $M\equiv-48,-44,15,-4,56,-55,-45,-61,26,49$ $\textup{(mod}$ $229\textup{)}$ and
$S_{r-1}^{(1)}\equiv-1\equiv S_{r-1}^{(9)}$ $\textup{(mod}$ $M\textup{)}$,
$S_{r-1}^{(2)}\equiv-8$ $\textup{(mod}$ $M\textup{)}$,
$S_{r-1}^{(3)}\equiv7$ $\textup{(mod}$ $M\textup{)}$,
$S_{r-1}^{(4)}\equiv21$ $\textup{(mod}$ $M\textup{)}$,
$S_{r-1}^{(5)}\equiv-15$ $\textup{(mod}$ $M\textup{)}$,
$S_{r-1}^{(6)}\equiv-20$ $\textup{(mod}$ $M\textup{)}$,
$S_{r-1}^{(7)}\equiv10$ $\textup{(mod}$ $M\textup{)}$ and
$S_{r-1}^{(8)}\equiv5$ $\textup{(mod}$ $M\textup{)}$;

$\mathrm{(ii)}$ $M\equiv-98,38,92,-69,112,-35,-77,-32$ $\textup{(mod}$ $229\textup{)}$ and
$S_{r-1}^{(1)}\equiv1\equiv S_{r-1}^{(9)}$ $\textup{(mod}$ $M\textup{)}$,
$S_{r-1}^{(2)}\equiv-8$ $\textup{(mod}$ $M\textup{)}$,
$S_{r-1}^{(3)}\equiv-7$ $\textup{(mod}$ $M\textup{)}$,
$S_{r-1}^{(4)}\equiv21$ $\textup{(mod}$ $M\textup{)}$,
$S_{r-1}^{(5)}\equiv15$ $\textup{(mod}$ $M\textup{)}$,
$S_{r-1}^{(6)}\equiv-20$ $\textup{(mod}$ $M\textup{)}$,
$S_{r-1}^{(7)}\equiv-10$ $\textup{(mod}$ $M\textup{)}$ and
$S_{r-1}^{(8)}\equiv5$ $\textup{(mod}$ $M\textup{)}$.\epro

\bpf{}\;We let $L=\mathbb{Q}(\zeta_{38})$ , then
$Norm_{L/\mathbb{Q}}(\pi)=(\pi\bar{\pi})^{\sum\limits_{i=1}^9\sigma_{2i-1}}=229$. Since
$M\equiv-48,-44,15,-4,56,-55,-45,-61,26,49,-98,38,92,-69,112,-35,-77,-32\pmod{229},\\
n>1$, we have
$\left(\frac{M}{\pi}\right)_{38}\equiv
M^{(229-1)/38}=M^6\equiv-4,16,-64,-26,42,-15,60,-68,\\
43,4,-42,15,-60,-44,-53,-17\equiv\zeta_{19},\;
\zeta_{19}^2,\;\zeta_{19}^3,\;\zeta_{19}^6,\;\zeta_{19}^8,\;\zeta_{19}^{10},\;
\zeta_{19}^{11},\;\zeta_{19}^{16},\;\zeta_{19}^{17},\;-\zeta_{19},\;-\zeta_{19}^8,\\
-\zeta_{19}^{10},
\;-\zeta_{19}^{11},\;-\zeta_{19}^{13},\;-\zeta_{19}^{14},\;-\zeta_{19}^{15}$
 $\textup{(mod}$ $\pi\textup{)}$, where $26^6\equiv49^6\equiv43$
$\textup{(mod}$ $229\textup{)}$ and $38^6\equiv92^6\equiv-42$ $\textup{(mod}$ $229\textup{)}$, thus we get
$\left(\frac{M}{\pi}\right)_{38}=\zeta_{19},\;
\zeta_{19}^2,\;\zeta_{19}^3,\;\zeta_{19}^6,\;\zeta_{19}^8,\;\zeta_{19}^{10},\;
\zeta_{19}^{11},\;\zeta_{19}^{16},\;\zeta_{19}^{17},\;-\zeta_{19},\;-\zeta_{19}^8,\\
-\zeta_{19}^{10},
\;-\zeta_{19}^{11},\;-\zeta_{19}^{13},\;-\zeta_{19}^{14},\;-\zeta_{19}^{15}$. Here we notice that
$F(x)=F_7(x)=x^9+x^8-8x^7-7x^6+21x^5+15x^4-20x^3-10x^2+5x+1$, which implies $a_1=-1$,
$a_2=-8$, $a_3=7$, $a_4=21$, $a_5=-15$, $a_6=-20$, $a_7=10$, $a_8=5$, $a_9=-1$.
 Hence by the necessity of Theorem \ref{maintheorem}, if
$M$ is prime and $M\equiv-48,-44,15,-4,56,-55,-45,-61,26,49$ $\textup{(mod}$ $229\textup{)}$,
we have
$S_{r-1}^{(1)}\equiv a_1=-1$ $\textup{(mod}$ $M\textup{)}$,
$S_{r-1}^{(2)}\equiv a_2=-8$ $\textup{(mod}$ $M\textup{)}$,
$S_{r-1}^{(3)}\equiv a_3=7$ $\textup{(mod}$ $M\textup{)}$,
$S_{r-1}^{(4)}\equiv a_4=21$ $\textup{(mod}$ $M\textup{)}$,
$S_{r-1}^{(5)}\equiv a_5=-15$ $\textup{(mod}$ $M\textup{)}$,
$S_{r-1}^{(6)}\equiv a_6=-20$ $\textup{(mod}$ $M\textup{)}$,
$S_{r-1}^{(7)}\equiv a_7=10$ $\textup{(mod}$ $M\textup{)}$,
$S_{r-1}^{(8)}\equiv a_8=5$ $\textup{(mod}$ $M\textup{)}$,
$S_{r-1}^{(9)}\equiv a_9=-1$ $\textup{(mod}$ $M\textup{)}$,
since every $\left(\frac{M}{\pi}\right)_{38}$ is a primitive $19$-th root of unity.
If $M$ is prime and  $M\equiv-98,38,92,-69,112,-35,-77,-32$ $\textup{(mod}$ $229\textup{)}$,
then we can obtain
$S_{r-1}^{(1)}\equiv1\equiv S_{r-1}^{(9)}$ $\textup{(mod}$ $M\textup{)}$,
$S_{r-1}^{(2)}\equiv-8$ $\textup{(mod}$ $M\textup{)}$,
$S_{r-1}^{(3)}\equiv-7$ $\textup{(mod}$ $M\textup{)}$,
$S_{r-1}^{(4)}\equiv21$ $\textup{(mod}$ $M\textup{)}$,
$S_{r-1}^{(5)}\equiv15$ $\textup{(mod}$ $M\textup{)}$,
$S_{r-1}^{(6)}\equiv-20$ $\textup{(mod}$ $M\textup{)}$,
$S_{r-1}^{(7)}\equiv-10$ $\textup{(mod}$ $M\textup{)}$,
$S_{r-1}^{(8)}\equiv5$ $\textup{(mod}$ $M\textup{)}$,
because $\left(\frac{M}{\pi}\right)_{38}$ is a primitive $38$-th root of unity.
This completes the proof of necessity.

Next by the sufficiency
of Theorem \ref{maintheorem} and the assumption, we have whatever $\mathrm{(i)}$ or
$\mathrm{(ii)}$ holds, $M$ is prime. This
completes the proof of sufficiency.
\qed\epf
%\smallskip
\textbf{Remark.}\quad $\mathrm{(i)}$ \;For $5\leq p\leq19$, every element $\alpha_i^{(k+1)}$, $1\leq i\leq(p-1)/2$ exists a polynomial recurrent relation about
the corresponding $\alpha_i^{(k)}$. By the definition of $S_k^{(j)}$, $1\leq j\leq(p-1)/2$,
we may get a recurrent relation for every $S_{k+1}^{(j)}$ on all elements $S_k^{(i)}$, $1\leq i\leq(p-1)/2$. Hence
there are $(p-1)/2$ many recurrent relations theoretically.

$\mathrm{(ii)}$ \;Once the recurrent relations are given, we will acquire the explicit primality tests for $p\geq5$.
 And we will show the definite recurrent relations of  $S_{k+1}^{(j)}$, $j=1,2$,
involved in Proposition \ref{2} for $p=5$ later. Then we could deduce the explicit primality test
for $p=5$. As to others $p\geq7$ the recurrent relations
may be very long and complicated, we don't write down them in the paper.

\section{Implementation and Computational results}\label{SecImple}
In this section we will verify the correctness of the algorithms related to Propositions \ref{1} and
\ref{2}. Let $G_n=6^{2^n}+1$ and $H_n=10^{2^n}+1$. First we make some preparations for the case $p=5$.
When $k\geq0$, the recurrent relations of $S_{k+1}^{(j)}$, $j=1,2$,
involved in Proposition \ref{2} can be obtained as follows.

 By the definition of $\alpha_1^{(k)}$ and $\alpha_2^{(k)}$, we have
  $$\alpha_1^{(k+1)}=(\alpha_1^{(k)})^{10}-10(\alpha_1^{(k)})^{8}
+35(\alpha_1^{(k)})^{6}-50(\alpha_1^{(k)})^{4}+25(\alpha_1^{(k)})^{2}-2.$$

Hence
$$\alpha_2^{(k+1)}=\sigma_3(\alpha_1^{(k+1)})=(\alpha_2^{(k)})^{10}-10(\alpha_2^{(k)})^{8}
+35(\alpha_2^{(k)})^{6}-50(\alpha_2^{(k)})^{4}+25(\alpha_2^{(k)})^{2}-2.$$

From the expression of $S_{k}^{(1)}$ and $S_{k}^{(2)}$ in Proposition \ref{2}, and under some computations, we get
\[
\begin{array}{ll}
S_{k+1}^{(1)}
&=(S_{k}^{(1)})^{10}-10(S_{k}^{(1)})^8S_{k}^{(2)}+35(S_{k}^{(1)})^6(S_{k}^{(2)})^2-50(S_{k}^{(1)})^4(S_{k}^{(2)})^3\\\\
&+25(S_{k}^{(1)})^2(S_{k}^{(2)})^4-10(S_{k}^{(1)})^8+80(S_{k}^{(1)})^6S_{k}^{(2)}-200(S_{k}^{(1)})^4(S_{k}^{(2)})^2\\\\
&+160(S_{k}^{(1)})^2(S_{k}^{(2)})^3-20(S_{k}^{(2)})^4+35(S_{k}^{(1)})^6-210(S_{k}^{(1)})^4S_{k}^{(2)}\\\\
&+315(S_{k}^{(1)})^2(S_{k}^{(2)})^2-70(S_{k}^{(2)})^3-50(S_{k}^{(1)})^4+200(S_{k}^{(1)})^2S_{k}^{(2)}\\\\
&-100(S_{k}^{(2)})^2-2(S_{k}^{(2)})^5+25(S_{k}^{(1)})^2-50S_{k}^{(2)}-4,
\end{array}
\]
and
\[
\begin{array}{ll}
S_{k+1}^{(2)}
&=(S_{k}^{(2)})^{10}+20(S_{k}^{(2)})^9-10(S_{k}^{(1)})^2(S_{k}^{(2)})^8+170(S_{k}^{(2)})^8-140(S_{k}^{(1)})^2(S_{k}^{(2)})^7\\\\
&+800(S_{k}^{(2)})^7+35(S_{k}^{(1)})^4(S_{k}^{(2)})^6-800(S_{k}^{(1)})^2(S_{k}^{(2)})^6+2275(S_{k}^{(2)})^6\\\\
&+300(S_{k}^{(1)})^4(S_{k}^{(2)})^5-2400(S_{k}^{(1)})^2(S_{k}^{(2)})^5+4004(S_{k}^{(2)})^5-50(S_{k}^{(1)})^6(S_{k}^{(2)})^4\\\\
&+1000(S_{k}^{(1)})^4(S_{k}^{(2)})^4-4050(S_{k}^{(1)})^2(S_{k}^{(2)})^4+4290(S_{k}^{(2)})^4-200(S_{k}^{(1)})^6(S_{k}^{(2)})^3\\\\
&+1600(S_{k}^{(1)})^4(S_{k}^{(2)})^3-3820(S_{k}^{(1)})^2(S_{k}^{(2)})^3+2640(S_{k}^{(1)})^3+25(S_{k}^{(1)})^8(S_{k}^{(2)})^2\\\\
&-320(S_{k}^{(1)})^6(S_{k}^{(2)})^2+1275(S_{k}^{(1)})^4(S_{k}^{(2)})^2-1880(S_{k}^{(1)})^2(S_{k}^{(2)})^2+825(S_{k}^{(2)})^2\\\\
&+20(S_{k}^{(1)})^8S_{k}^{(2)}-160(S_{k}^{(1)})^6S_{k}^{(2)}+420(S_{k}^{(1)})^4S_{k}^{(2)}-400(S_{k}^{(1)})^2S_{k}^{(2)}-2(S_{k}^{(1)})^{10}\\\\
&+20(S_{k}^{(1)})^8-70(S_{k}^{(1)})^6+100(S_{k}^{(1)})^4-50(S_{k}^{(1)})^2+100S_{k}^{(2)}+4.
\end{array}
\]
With the above two recurrent formulas, we can easily reach the explicit primality test for numbers $H_n$.

\textbf{Remark.} Though the expression of the above two recurrent relations is a bit long. Its corresponding algorithm has
the same complexity of running time compared with the primality test of $H_n$ in \cite{wi}, which runs in polynomial time in log$_{2}(M)$.

 We implement our two algorithms of $p=3$ and $p=5$ in Magma \cite{magma}. And our program is run on a personal computer with Intel Core i5-3470 3.20GHz CPU and 4GB memory.

 We verified the correctness of our program by comparing with the results in \cite{ri2}. Since the growth of numbers $G_n$ and $H_n$ is
very fast on index $n$. When $n\geq15$ the computation involved in our program is out of memory in our personal computer. If we deal with
better and more efficient representation of lager integers, we may test the primality of bigger $G_n$ or $H_n$. But
this is not the stress of our paper. And we will not mention it any more here. We verified all numbers for $p=3, 5$ in
the range $1\leq n<15$ and found no mistakes. Note that the assumption about equation  $x^4\equiv1$
$\textup{(mod}$ $5^r\textup{)}$ in Proposition \ref{2} holds for $H_n$, $1\leq n<15$, by applying the algorithm in \cite{dl}.
 The prime numbers are rare on such $G_n$ and $H_n$. We list the following two tables to show all
the cases and their corresponding cost time.

\newpage

\begin{table}
\caption{Primality of $G_n=6^{2^n}+1(p=3)$}
\begin{tabular}{|c|c|c|l|}
\hline n     & $G_n$     & Primality & Time(s)\\
\hline
1     & 37    & yes       & 0.011       \\
\hline
2     & 1297  & yes       & 0.015      \\
\hline
3 to 10 & -     & no        & 0.921        \\
\hline
11    & -     & no        & 3.931  \\
\hline
12    & -     & no        & 23.228 \\
\hline
13    & -     & no        & 139.293\\
\hline
14    & -     & no        & 738.805\\
\hline
\end{tabular}
\end{table}

\begin{table}
\caption{Primality of $H_n=10^{2^n}+1(p=5)$}
\begin{tabular}{|c|c|c|l|}
\hline n     & $H_n$     & Primality & Time(s)\\
\hline
1     & 101    & yes       &  0.015      \\
\hline
2 to 10 & -     & no        & 7.909       \\
\hline
11    & -     & no        & 37.004  \\
\hline
12    & -     & no        & 204.579 \\
\hline
13    & -     & no        & 1180.226 \\
\hline
14    & -     & no        & 6576.924  \\
\hline
\end{tabular}
\end{table}

%\vspace{0.4cm}
\vspace{0.6cm}
\centerline{\bf Acknowledgments}
\vspace {0.6cm} The authors are very grateful to Yupeng Jiang and Chang Lv for their useful
discussion. The work of this paper
was supported by the NNSF of China (Grants Nos. 11071285,
61121062), 973 Project (2011CB302401) and the National Center for Mathematics and Interdisciplinary Sciences, CAS.

%\vskip0.2in\noindent{\bf References} \vskip0.1in

\end{document}